\documentclass[11pt,amsthm,reqno]{amsart}
\usepackage{amsmath,amssymb}
\usepackage{amscd}

\newcommand{\C}{{\mathbb C}}
\newcommand{\Z}{{\mathbb Z}}
\newcommand{\Q}{{\mathbb Q}}
\newcommand{\A}{{\mathbb A}}

\newcommand{\Spec}{{\rm Spec}\;}
\newcommand{\Ker}{{\rm Ker}\:}

\newcommand{\codim}{{\rm codim}}

\newcommand{\Sing}{{\rm Sing}}
\newcommand{\pl}{{\rm pl}}
\newtheorem{thm}{Theorem}[section]
\newtheorem{cor}[thm]{Corollary}
\newtheorem{lem}[thm]{Lemma}
\newtheorem{prop}[thm]{Proposition}
\newcommand{\Proof}{\noindent{\bf Proof.}\quad}
\newcommand{\Remark}{\noindent{\bf Remark.}\quad}
\newcommand{\Remarks}{\noindent{\bf Remarks.}\quad}

\begin{document}
\title[Faithfully flat quotient morphisms by $G_a$-actions]
{Faithfully flat quotient morphisms by $G_a$-actions on factorial affine varieties}
\author{Kayo Masuda}
\address{Department of Mathematical Sciences, 
School of Science, Kwansei Gakuin University \\
1 Gakuen Uegahara, Sanda 669-1330, Japan}
\email{kayo@kwansei.ac.jp}
\thanks{Supported by JSPS KAKENHI Grant number JP20K03570, JP25K06945.}
\date{}

\begin{abstract}
Let $X$ be a factorial complex affine variety of dimension $\ge 3$ 
with an algebraic action of the additive group $G_a$. 
Let $\pi : X \to Y$ be the algebraic quotient morphism where we assume $Y$ is an affine variety. 
When $\pi$ is faithfully flat, we investigate $\pi$ by $G_a$-equivariant affine modifications and 
give criteria for $\pi$ to be a trivial $\A^1$-bundle.
For a smooth acyclic fourfold $X$ with a free $G_a$-action and
a $G_a$-equivariant $\A^3$-fibration $f : X \to \A^1$ where $G_a$ acts trivially on $\A^1$, 
we give a criterion for the algebraic quotient $Y$
to be isomorphic to $\A^3$ with $f$ as a coordinate. 
Together with a criterion for $\pi : X \to Y$ to be a trivial $\A^1$-bundle,
we obtain a sufficient condition for $X\cong Y\times \A^1\cong \A^4$. 
\end{abstract}

\maketitle

\section{Introduction} 

Let $X=\Spec B$ be an affine algebraic variety defined over $\C$
with a non-trivial action of the additive group $G_a$. 
There is a non-trivial locally nilpotent derivation (abbreviated to lnd) $\delta$ on $B$
associated with the $G_a$-action on $X$. 
Let $A=\Ker \delta$. 
The subalgebra $A$ of $B$ is finitely generated over $\C$ if $\dim X\le 3$
by the Zariski finiteness theorem \cite{Zar},
however, if $\dim X >3$, $A$ is not necessarily an affine $\C$-domain. 
The morphism $\pi : X \to Y=\Spec A$ of affine schemes induced by
the inclusion $A \hookrightarrow B$ is called the algebraic quotient morphism,
which we call simply the quotient morphism in this article. 
The quotient morphism $\pi : X \to Y$ is an $\A^1$-fibration 
whose general closed fiber is a $G_a$-orbit.
The quotient morphism $\pi$ is not necessarily equi-dimensional nor surjective 
even though $Y$ is an affine variety and 
the $G_a$-action on $X$ is free, i.e., $X$ has no $G_a$-fixed points 
(for example, see Winkelmann \cite{Win}). 

An affine $n$-space $\A^n$ for $n\le 3$ can be characterized by a non-trivial $G_a$-action. 
If an affine curve $X$ admits a $G_a$-action, $X$ is isomorphic to affine line $\A^1$. 
If an affine surface $X=\Spec B$ is factorial, i.e., $B$ is a UFD,
the group $B^*$ of invertible elements of $B$ is $\C^*$ and
$X$ admits a $G_a$-action, then $X\cong \A^2$ by Miyanishi \cite{Miy3}. 
If an affine threefold $X=\Spec B$ is smooth, acyclic, i.e., $H_i(X,\Z)=0$ for every $i>0$, and 
there exists an $\A^2$-fibration $f : X \to \A^1$ such that every closed fiber is factorial,
then $X\cong \A^3$ (Miyanishi \cite{Miy4}, Kaliman \cite{Kal2}, cf. Theorem \ref{Theorem 4.1}). 
Note that any smooth acyclic affine variety $Z=\Spec R$ is factorial and $R^*=\C^*$ 
by Fujita \cite[1.18-1.20]{Fuj} (see also \cite[3.2]{Kal}). 
In the case that $X=\A^n$ for $n\le 3$, $Y=\A^{n-1}$ and
the quotient morphism $\pi : \A^n \to \A^{n-1}$ is faithfully flat, hence
equi-dimensional and surjective (see \cite{Bon}, \cite{DK} for $n=3$). 
Further, every free $G_a$-action on $\A^n$ for $n\le 3$ is a translation, i.e., 
$\pi : X \to Y$ is a trivial $\A^1$-bundle (for $n=3$, see \cite{Kal1}). 
For $n=4$, a proper $G_a$-action on $\A^4$ preserving a coordinate is a translation
by Kaliman \cite{Kal3}.
A $G_a$-action on $X$ is proper
iff the morphism $G_a\times X\ni (t,x)\mapsto (x, t\cdot x)\in X\times X$ is proper. 
A proper $G_a$-action is free. 
However, for $n\ge 4$, the situation is completely different.  
In fact, the quotient morphism by a proper $G_a$-action is not necessarily
equi-dimensional nor surjective as shown by one of Winkelmann's examples \cite{Win}. 
The example is a $G_a$-action on $\A^5$ preserving two coordinates which is proper but not 
a translation. 
Hence if we investigate
when the quotient morphism $\pi : \A^n \to Y$ is a trivial $\A^1$-bundle for $n\ge 4$,
we should consider under some condition, say, $\pi$ is faithfully flat or, 
equi-dimensional and surjective.  

In this article, we study the quotient morphism $\pi : X \to Y$ by a $G_a$-action 
on a factorial affine variety $X$ of dimension $\ge 3$. 
We investigate $\pi$ under the assumption that $Y$ is an affine variety and 
$\pi$ is faithfully flat. 
When $\pi$ is faithfully flat, 
the set $\Sing (\pi)=\{Q \in Y \mid \pi^*(Q)\not\cong \A^1\}$ is empty or 
of pure codimension one in $Y$ (\cite{Ma2}, Proposition \ref{Proposition 2.5} and its remark). 
If $\Sing (\pi)=\emptyset$, $\pi$ is a trivial $\A^1$-bundle ({\it ibid}).
In the case that $X$ is smooth and the $G_a$-action on $X$ is free,
if $\pi : X \to Y$ is equi-dimensional and surjective, then $\pi$ is faithfully flat
(Proposition \ref{Proposition 2.5}).
We give criteria for a faithfully flat $\pi : X \to Y$ to be a trivial $\A^1$-bundle
in view of the geometry of $G_a$-actions (Theorems \ref{Theorem 3.5}, \ref{Theorem 3.6}). 
We also give a simple proof to the characterization of $\A^3$ stated above
(Theorem \ref{Theorem 4.1}). 
Let $X$ be a smooth acyclic fourfold with a free $G_a$-action and
a $G_a$-equivariant $\A^3$-fibration $f: X \to \A^1$ where $G_a$ acts on $\A^1$ trivially. 
Using the characterization of $\A^3$, we show that $Y\cong \A^3$ and
$f$ is a coordinate of $Y=\A^3$ under the assumptions that 
the quotient morphism $\pi : X \to Y$ is equi-dimensional and surjective,
$Y$ is an acyclic affine variety and every closed fiber of the morphism $f_Y : Y \to \A^1$
induced by $f$ is factorial (Theorem \ref{Theorem 4.3}). 
Further, if a criterion for $\pi$ to be a trivial $\A^1$-bundle is satisfied, 
then $X\cong Y\times \A^1\cong \A^4$.  

The construction of this article is as follows. 
After reviewing quotient morphisms by $G_a$-actions in Section 2,
we analyze faithfully flat $\pi : X \to Y$ by $G_a$-equivariant affine modifications 
and describe the $G_a$-action on $X$ in Lemmas \ref{Lemma 3.3} and \ref{Lemma 3.4}.
The criteria for $\pi$ to be a trivial $\A^1$-bundle are given 
in Theorems \ref{Theorem 3.5} and \ref{Theorem 3.6}. 
In Section 4, we consider the case that $X$ is smooth and acyclic. 
We give a simple proof to the characterization of $\A^3$ in Theorem \ref{Theorem 4.1}.
Finally, we study a smooth acyclic fourfold $X$ with a free $G_a$-action and
a $G_a$-equivariant $\A^3$-fibration $f: X \to \A^1$. 
We give a criterion for $Y\cong \A^3$ in Theorem \ref{Theorem 4.3}. 
By a criterion for $\pi : X \to Y$ to be a trivial $\A^1$-bundle, 
we obtain a sufficient condition for $X\cong Y\times \A^1\cong \A^4$. 

\medskip
\bigskip

\section{Quotient morphisms by $G_a$-actions}

We review some facts on $G_a$-actions and the quotient morphisms. 
For general references, we refer to \cite{GMM3} and \cite{Fre}. 

Let $X=\Spec B$ be an affine variety with a non-trivial $G_a$-action
associated with an lnd $\delta$ on the affine $\C$-domain $B$. 
Let $A=\Ker \delta$ and $Y=\Spec A$.
The subalgebra $A$ of $B$ is factorially closed in $B$, i.e.,
for $b_1, b_2 \in B$, $b_1b_2 \in A\setminus \{0\}$ implies $b_1, b_2 \in A$.
Hence if $X$ is factorial (resp. normal), then $Y$ is factorial (resp. normal) as well. 
Further, $A^*=B^*$. 
There exists $z\in B$, which is transcendental over $A$, 
such that $\delta(z)=a\in A\setminus \{0\}$ and $B[a^{-1}]=A[a^{-1}][z]$. 
Such an element $z\in B$ is called a local slice of $\delta$. 
If $\delta(z)=1$, $z$ is called a slice. 
The ideal $\delta(B)\cap A$ of $A$ is called the plinth ideal and denoted by $\pl (\delta)$.  
The quotient morphism $\pi : X \to Y$ is defined by the inclusion $A \hookrightarrow B$.  
Note that the closed subset $V(\pl(\delta))$ of $Y$ contains $\Sing (\pi)$. 

An lnd $\delta$ on $B$ is called {\it irreducible}
iff $\delta(B)\subset bB$ for $b \in B$ implies $b \in B^*$. 
There exist an irreducible lnd $\tilde \delta$ on $B$ and 
an element $a\in \Ker \tilde\delta$ such that $\delta=a\tilde\delta$ 
\cite[Principle 7, Proposition 2.2]{Fre}. 
Note that $\Ker \delta=\Ker \tilde \delta$. 
Hence the quotient morphism induced by the lnd $\tilde \delta$ coincides with the one by $\delta$. 

The fixed point locus $X^{G_a}$ is defined by the ideal generated by $\delta(B)$. 
Hence if $B$ is factorial and $\delta$ is irreducible, then $X^{G_a}$ has codimension $>1$. 
By a result of Bia{\l}ynicki-Birula \cite[Corollary 4]{BB}, $X$ has no isolated fixed points. 

\begin{lem}\label{Lemma 2.1}
Let $X$ be a factorial affine surface with a $G_a$-action associated with an irreducible lnd. 
Then $X$ is smooth and has no fixed points.
\end{lem}
\Proof
Since $\dim X^{G_a}<1$ as remarked above, it follows from \cite{BB} that $X^{G_a}=\emptyset$. 
The $G_a$-action $\sigma : G_a \times X \to X$ induces 
the automorphism $X\ni P \mapsto \sigma(t,P) \in X$ for $t\in G_a$. 
If $P\in X$ is a singular point, then $\sigma(t,P)$ is a singular point as well.
Since $X$ is normal, the singular locus of $X$ is of codimension $>1$. 
Hence if a singular point of $X$ exists, it must be a fixed point. 
It follows from $X^{G_a}=\emptyset$ that $X$ is nonsingular.
\qed

\medskip

Let $q : X \to Y$ be an $\A^1$-fibration.
A fiber $q^*(Q)$ over a point $Q\in Y$ is {\it singular} if it is not isomorphic to $\A^1$. 
Let $F$ be a singular fiber of $q$ and $F=\sum_{i=1}^rF_i$ be the decomposition into
irreducible components. 
An irreducible component $F_i$ is a {\it multiple} fiber component of $q$ 
if the multiplicity of $F_i$ is greater than one.

\begin{thm}{\rm (\cite[Theorem 4.11]{GKMMR}, \cite[Corollary 3.1.39]{GMM3})}\label{Theorem 2.2}
Let $X$ be a smooth affine variety with a non-trivial $G_a$-action 
and let $\pi: X \to Y$ be the quotient morphism where we assume $Y$ is an affine variety. 
Suppose that $\pi$ is equi-dimensional. 
Then the following assertions hold. 

\begin{enumerate}
\item[(1)] Every multiple fiber component of $\pi$ is contained in $X^{G_a}$.
\item[(2)] Every reduced closed fiber $F\neq \emptyset$ of $\pi$ consists of
a disjoint union of $\A^1$'s. 
\end{enumerate} 
\end{thm}
\Proof
(1) follows from \cite[Theorem 4.11]{GKMMR}.

(2) Since $\pi$ is equi-dimensional, every closed fiber of $\pi$ has pure dimension one.
Then the assertion follows from \cite[Corollary 3.1.39]{GMM3}.
\qed

\medskip

If the $G_a$-action on $X$ is free in Theorem \ref{Theorem 2.2}, 
then there are no multiple fiber components of $\pi$ and every closed fiber of $\pi$ is reduced.   
Hence every closed fiber consists of a disjoint union of $\A^1$'s. 

If $X$ is a normal affine surface, then $\pi$ is surjective by \cite{GMM}. 
If $X$ is a factorial affine $n$-fold for $n \le 3$, then $\pi$ is equi-dimensional ({\it ibid.}). 
Hence if $X$ is a smooth acyclic affine threefold, $\pi$ is equi-dimensional,
and also surjective by Kaliman \cite[Remark 3.3]{Kal1}. 
Kaliman and Saveliev \cite{KS} showed that if $X$ is a smooth contractible affine threefold
with a free $G_a$-action, 
then $\pi :X \to Y$ is a trivial $\A^1$-bundle and $X\cong Y\times \A^1$.  
As remarked in Section one, when $\dim X>3$, 
$\pi$ is not necessarily equi-dimensional nor surjective even if the action is free. 

As for a faithfully flat $\A^1$-fibration, the following holds.

\begin{prop}{\rm (\cite[Corollary 2.1.5]{GMM3})}\label{Proposition 2.3}
Let $q : X \to Y$ be a faithfully flat fibration of affine varieties.
If $Y$ is normal and every closed fiber $q^*(Q)$ for $Q\in Y$ is isomorphic to $\A^1$, 
then $q$ is an $\A^1$-bundle.
\end{prop}

As for the smoothness of the base space $Y$ of the quotient morphism $\pi : X \to Y$,
we have the following result.
The proof is essentially the same as the one of \cite[Lemma 2.22]{GKMMR2}. 

\begin{prop}{\rm (cf. \cite[Lemma 2.22]{GKMMR2})}\label{Proposition 2.4}
Let $q : X \to Y$ be a dominant morphism from an algebraic variety $X$ of dimension $n+r$ to
an algebraic variety $Y$ of dimension $n$. 
Let $Q\in Y$.
Suppose that $Y$ is normal at $Q$ and the fiber $q^*(Q)$ has an irreducible component $F$
which is reduced and $r$-dimensional. 
If there exists a point $P\in F$ at which $X$ and $F$ are smooth, 
then $Q$ is a smooth point of $Y$.  
\end{prop}

\medskip

By Theorem \ref{Theorem 2.2} and Proposition \ref{Proposition 2.4}, the following holds.

\begin{prop}\label{Proposition 2.5}
Let $X$ be a smooth affine variety of dimension $\ge 3$ with a non-trivial $G_a$-action and 
let $\pi :X \to Y$ be the quotient morphism where we assume $Y$ is an affine variety. 
Suppose that $\pi$ is surjective and equi-dimensional. Then the following assertions hold. 
\begin{enumerate}
\item[(1)] If the $G_a$-action on $X$ is free, then $Y$ is smooth and $\pi$ is faithfully flat. 
\item[(2)] If $Y$ is smooth and $\Sing (\pi)=\emptyset$, $\pi$ is an $\A^1$-bundle. 
\item[(3)] Suppose that $X$ is factorial and $Y$ is smooth.
Then the plinth ideal $\pl (\delta)$ is principal.   
If $\Sing (\pi)=\emptyset$, then $\pi$ is a trivial $\A^1$-bundle.   
If $\Sing (\pi)\neq \emptyset$, then the closure $\overline{\Sing (\pi)}$ coincides with
$V(\pl (\delta))$ and $\Sing (\pi)$ is of pure codimension one in $Y$. 
\end{enumerate}
\end{prop}
\Proof
(1) If the $G_a$-action is free, then every closed fiber is reduced by Theorem \ref{Theorem 2.2}. 
Hence $Y$ is smooth by Proposition \ref{Proposition 2.4}.
Then $\pi$ is faithfully flat \cite{Matsu}. 

(2) Since $Y$ is smooth, $\pi$ is faithfully flat ({\it ibid.}). 
The assertion follows from Proposition \ref{Proposition 2.3}. 

(3) Let $\delta$ be the lnd associated with the $G_a$-action. 
Since $Y$ is smooth, $\pi$ is faithfully flat. 
Then $\pl (\delta)$ is principal by \cite{Ma2}.
If $\Sing (\pi)=\emptyset$, then $\pi$ is a trivial $\A^1$-bundle by \cite[Corollary 3.4]{Ma2}.
Suppose $\Sing (\pi)\neq \emptyset$. 
We may assume $\delta$ is irreducible. 
Since $\Sing (\pi) \subset V(\pl(\delta))$ and
each irreducible component of $V(\pl(\delta))$ is contained in $\overline{\Sing (\pi)}$
by \cite[Theorem 3.3]{Ma2},
we have $\overline{\Sing (\pi)}=V(\pl (\delta))$ and the assertion follows. 
\qed

\medskip

\Remark Let $X=\Spec B$ be a factorial affine variety with a $G_a$-action
associated with an irreducible lnd $\delta$ and let $A=\Ker \delta$.
If $B$ is faithfully flat over $A$, then the plinth ideal $\pl(\delta)$ is principal
(cf. \cite[Lemma 2.1]{Ma2}).
Then $\overline{\Sing (\pi)}=V(\pl(\delta))$
by the same argument as in Proposition \ref{Proposition 2.5} (3). 

\bigskip

By the same argument as in the proof of \cite[Theorem 2.1]{GMM2},
we have the following result.

\begin{prop}{\rm (cf.\cite[Theorem 2.1]{GMM2})}\label{Proposition 2.6}
Let $X$ be a smooth acyclic $n$-fold for $n\ge 3$ and 
let $p : X \to \A^1$ be a faithfully flat morphism. 
Suppose that $p$ restricts to a trivial $\A^{n-1}$-bundle $p|_{p^{-1}(C)} : p^{-1}(C) \to C$ 
where $C=\A^1 \setminus \{P_1,\ldots, P_r\}$ is an affine line $\A^1$
with $r$ points deleted off and that $p^{-1}(P_i)$ is a smooth (irreducible) variety for
$1 \le i \le r$. 
Then $p^{-1}(P_i)$ is acyclic for $1 \le i \le r$.
\end{prop}
\Proof
Let $X'=\bigcup_{i=1}^r p^{-1}(P_i)$.  
We have a long exact sequence of integral cohomology groups with compact support for a pair
$(X, X')$
\begin{align*}
\cdots \to H^{2i}_c(X) \to & H^{2i}_c(X')\to H^{2i+1}_c(X,X')\to \\
& H^{2i+1}_c(X)\to H^{2i+1}_c(X')\to H^{2i+2}_c(X,X')\to \cdots .
\end{align*}
By the Poincar\'e duality, we have $H^{2i}_c(X)\cong H_{2n-2i}(X)=0$ and
$H^{2i+1}_c(X)\cong H_{2n-2i-1}(X)=0$ for $i \le n-1$.
By the Lefschetz duality, $H^{2i+1}_c(X,X')\cong H_{2n-2i-1}(X\setminus X')\cong H_{2n-2i-1}(C)=0$ 
and $H^{2i+2}_c(X,X')\cong H_{2n-2i-2}(X\setminus X')\cong H_{2n-2i-2}(C)=0$ for $i \le n-2$.
Hence $H^{2i}_c(X')\cong H_{2n-2-2i}(X')=0$ and $H^{2i+1}_c(X')\cong H_{2n-2i-3}(X')=0$
for $i\le n-2$. Thus $H_{j}(X')=0$ for $j\ge 1$, and the assertion follows.
\qed

\medskip

\Remark Let $X$ be a smooth affine $n$-fold for $n\ge 3$ and
let $p : X \to \A^1$ and $P_1, \ldots, P_r$ be as in Proposition \ref{Proposition 2.6}.
Then by the same argument as in the proof of \cite[Theorem 2.3]{GMM2},
$X$ is acyclic if $p^*(P_i)$ is reduced and $p^{-1}(P_i)$ is acyclic for every $i$. 
Further, since $\pi_1(X)$ is trivial by Nori \cite{Nor},
$X$ is contractible by a theorem of J. H. C. Whitehead. 

\medskip
\bigskip

\section{$G_a$-equivariant affine modifications}

We investigate $G_a$-actions on a factorial affine variety
by $G_a$-equivariant affine modifications. 
We begin by recalling the definition of affine modification.
For a general reference, we refer to \cite{KZ}. 

Let $B$ be an affine $\C$-domain.
Let $R$ be a subdomain of $B$, $I$ an ideal of $R$, and $f$ a nonzero element of $I$. 
The subalgebra of the quotient field $Q(R)$ generated by the elements $a/f$ for $a\in I$
is called the {\it affine modification} of $R$ along $f$ with center $I$ and
denoted by $R[f^{-1}I]$. 
Suppose that $B$ is equipped with an lnd $\delta$ restricting to $R$ such that $\delta(f)=0$. 
Then the ideal $J=R\cap fB$ of $R$ is $\delta$-stable, i.e., $\delta(J)\subset J$, 
and $\delta$ restricts to $R[f^{-1}J]\subset B$. 
The inclusions $R \hookrightarrow R[f^{-1}J] \hookrightarrow B$ yield 
a sequence of $G_a$-equivariant birational morphisms 
$\Spec B \to \Spec R[f^{-1}J] \xrightarrow{\sigma} \Spec R$.
The reduced divisor $(f\circ \sigma)^*(0)$ is called the exceptional divisor of $\sigma$. 

In the following, $B$ is a factorial affine $\C$-domain
equipped with an irreducible lnd $\delta$ such that $A=\Ker \delta$ is an affine $\C$-domain. 
Suppose that 
\begin{equation}
\pl(\delta)=f^p\beta A \label{(1)}
\end{equation}  
where $f$ is a prime element of $A$, $p>0$, and $\beta \in A$ and $f$ are coprime. 
Let $z \in B$ be a local slice such that $\delta(z)=f^p\beta$. 
Then $z$ is not divisible by $f$ nor any prime factor of $\beta$. 
Since $B[f^{-1}\beta^{-1}]=A[f^{-1}\beta^{-1}][z]$ and $B$ is finitely generated over $\C$,
there exist positive integers $\mu$ and $\gamma$
such that
\begin{equation}
B=A[z][f^{-\mu}\beta^{-\gamma}\tilde J] \label{(2)}
\end{equation}  
where $\tilde J=A[z]\cap f^\mu\beta^\gamma B$.
Let
$$I=A[z]\cap f^{\mu}B \quad \text{and} \quad B^{(1)}=A[z][f^{-\mu}I].$$
Then
$$B=B^{(1)}[\beta^{-\gamma}J] \quad \text{where} \quad J=B^{(1)}\cap \beta^\gamma B.$$
Since $I$, $J$ and $\tilde J$ are $\delta$-stable ideals,
we obtain a sequence of $G_a$-equivariant affine modifications
\begin{equation}
B_0:=A[z]\subset B^{(1)}\subset B. \label{(3)}
\end{equation}  
Note that $B_0^{G_a}=(B^{(1)})^{G_a}=B^{G_a}=A$. 
Let $X=\Spec B$, $X^{(1)}=\Spec B^{(1)}$ and $X_0=\Spec B_0=Y\times \A^1$ where $Y=\Spec A$.
The sequence (\ref{(3)}) yields a sequence of $G_a$-equivariant birational morphisms over $Y$
\begin{equation*}
X \xrightarrow{\tau} X^{(1)} \xrightarrow{\sigma} X_0=Y\times \A^1 
\end{equation*}  
where $\sigma$ is the $G_a$-equivariant birational morphism induced by
the inclusion $B_0 \hookrightarrow B^{(1)}$. 
Note that $B[\beta^{-1}]=B^{(1)}[\beta^{-1}]$. 

We shall investigate the case that $\beta \in A^*$.  
Then $B=B^{(1)}$ and $X=X^{(1)}$. We may assume $\delta(z)=f^p$. 
Note that $f$ is a prime element of $B$ since $A$ is factorially closed in $B$. 
The exceptional divisor of $\sigma$ is $E=\Spec B/fB$. 
Let $\pi : X \to Y$ be the quotient morphism. 
Then $\pi (E)=\Spec \overline A$ where $\overline A=A/fA$.   

For $i\ge 1$, let $I_{0,i}=B_0\cap f^i B$.
Then $I_{0,i}$ is a $\delta$-stable ideal of $B_0=A[z]$ for every $i$ and 
$I_{0,1}=B_0\cap fB$ is a prime ideal.
The residue ring $B_0/fB_0$ is identified with a polynomial ring $\overline A[\overline z]$
over $\overline A$ where $\overline z$ is the residue class of $z$. 
Let $\overline I_{0,1}$ be the image of $I_{0,1}$ by the surjection 
$B_0\to B_0/fB_0=\overline A[\overline z]$. Then $\overline I_{0,1}\neq (0)$ by \cite{Ma2}.
Let $\overline B_0=B_0/I_{0,1}$. Note that
$$\overline A \subset \overline B_0 \subset B/fB$$
and $\sigma(E)=\Spec \overline B_0$. 
We have the following. 

\begin{lem}{\rm (cf. \cite[Lemma 4.2]{Ma2})}\label{Lemma 3.1}
Let $g\in I_{0,1}$ be an element which maps to a generator of 
the principal ideal $\overline I_{0,1}\otimes_{\overline A} K$ of $K[\overline z]$ by the map
$B_0\to B_0/fB_0=\overline A[\overline z]\hookrightarrow K[\overline z]$
where $K=Q(\overline A)$ is the quotient field of $\overline A$. 
Write $g\in B_0$ as 
$$g=f^{\ell_0}y_1$$
where $\ell_0\ge 1$ and $y_1\in B\setminus fB$.
Let $B_1=B_0[f^{-\ell_0}I_{0,\ell_0}]$. 
\begin{enumerate}
\item[(1)] There exists $c_0 \in A \setminus fA$ such that the ideal 
$\overline I_{0,1}\otimes_{\overline A}\overline A[\overline c_0^{-1}]$ of 
$\overline A[\overline c_0^{-1}][\overline z]$ is generated by $\overline g$ and
$I_{0,1}\otimes_A A[c_0^{-1}]=(f, g)B_0[c_0^{-1}]$. 
We have
$$B_1[c_0^{-1}]=B_0[c_0^{-1}][y_1]=A[c_0^{-1}][z,y_1]$$ and 
$$(B_1/fB_1)[\overline c_0^{-1}]=\overline B_0[\overline c_0^{-1}][\overline y_1]
=(\overline B_0[\overline c_0^{-1}])^{[1]}.$$ 
\item[(2)] $\deg_{\overline z}\overline g>1$ and $\overline{g'(z)}\notin \overline I_{0,1}$. 
\item[(3)] $[Q(\overline B_0):Q(\overline A)]>1$. 
\item[(4)] We have $p\ge \ell_0$ and  
$$\delta(y_1)=f^{p-\ell_0}g'(z).$$
Hence $\delta(B_1)\subset f^{p-\ell_0}B$.   
\item[(5)] Suppose that $\overline A$ is factorial. 
Take $g\in I_{0,1}$ so that $\overline g\in B_0/fB_0=\overline A[\overline z]$ 
is primitive over $\overline A$. Then $I_{0,1}=(f, g)B_0$ and $B_1=B_0[y_1]$ is factorial. 
Hence $B_1/fB_1=\overline B_0[\overline y_1]={\overline B_0}^{[1]}$ and
$$\overline B_0=B_0/I_{0,1}=\overline A[\overline z]/\overline g \overline A[\overline z].$$
\end{enumerate}
\end{lem}
\Proof 
(1) The element $\overline g\in \overline A[\overline z]$ has the minimal degree
with respect to $\overline z$ in $\overline I_{0,1}$.
Note that $\overline g\not\in \overline A$. 
Let $I_{0,1}=(f, h_1, h_2, \ldots,h_r)B_0$ for $h_1, \ldots, h_r \in B_0$. 
We may assume $h_1=g$ and $\{f, h_1=g, h_2, \ldots, h_r\}$ is a minimal set of generators.
For $2 \le i \le r$, we have $\overline a_i \overline h_i=\overline g\overline q_i$ 
in $B_0/fB_0=\overline A[\overline z]$ for $a_i \in A\setminus fA$ and $q_i \in B_0$. 
Let $c_0=a_2\cdots a_r \in A\setminus fA$. Then 
$\overline I_{0,1}\otimes_{\overline A}\overline A[\overline c_0^{-1}]$ is generated
by $\overline g$ and we have $I_{0,1}\otimes_A A[c_0^{-1}]=(f, g)B_0[c_0^{-1}]$.
By \cite{Ma2}, we have $I_{0,\ell_0}\otimes_A A[c_0^{-1}]=(f^{\ell_0},g)B_0[c_0^{-1}]$. 
Hence $B_1[c_0^{-1}]=B_0[c_0^{-1}][y_1]$. 
The assertion on $(B_1/fB_1)[\overline c_0^{-1}]$ follows from \cite[Lemma 4.2]{Ma2}. 

(2) The first assertion follows from \cite[Lemma 4.2]{Ma2}. 
As for the second assertion, 
note that $\overline{g'(z)}={\overline g}'(\overline z)$ is a nonzero element of
$\overline A[\overline z]$ and 
$\deg_{\overline z}{\overline g}'(\overline z)<\deg_{\overline z}\overline g(\overline z)$. 
Since $\overline g\in \overline A[\overline z]$ has the minimal degree in $\overline I_{0,1}$,
$\overline{g'(z)}$ is not contained in $\overline I_{0,1}$. 

(3) Note that $\overline B_0\otimes_{\overline A}K\cong K[\overline z]/\overline g K[\overline z]$ and $Q(\overline B_0)=K(\overline z)$.
The polynomial $\overline g\in \overline A[\overline z]\subset K[\overline z]$ is
the minimal polynomial of
$\overline z\in \overline B_0\subset Q(\overline B_0)=K(\overline z)$ over $K$ up to $K^*$. 
Hence by (2)
$$[Q(\overline B_0):Q(\overline A)]=[K(\overline z): K]=\deg_{\overline z}\overline g>1.$$

(4) By $f^{\ell_0} y_1=g \in A[z]$, we have
$f^{\ell_0}\delta(y_1)=g'(z) \cdot \delta(z)=f^p g'(z)$.
If $\ell_0>p$, then $g'(z)\in I_{0,1}$ and 
hence $\overline{g'(z)}\in \overline I_{0,1}$, which is a contradiction. 
Thus $p\ge \ell_0$ and we obtain the expression of $\delta(y_1)$.
Let $b_1\in B_1\subset B_1[c_0^{-1}]=A[c_0^{-1}][z,y_1]$. 
Since $\delta(b_1)\in f^{p-\ell_0}A[c_0^{-1}][z,y_1]$, the assertion follows. 

(5) Note that $\overline g$ is irreducible in $B_0/fB_0=\overline A[\overline z]$. 
In the argument in (1), we have $\overline h_i=\overline g \cdot \overline \xi_i$ 
for $\overline \xi_i\in \overline A[\overline z]$ since $\overline g$ is irreducible. 
Thus $\overline I_{0,1}$ is generated by $\overline g$ and we have $I_{0,1}=(f,g)B_0$. 
By \cite{Ma2}, $B_1=B_0[y_1]\cong B_0[Y]/(f^{\ell_0}Y-g)$ 
where $Y$ is an indeterminate. 
The factoriality of $B_1$ follows from Nagata's criterion for factoriality \cite{Nag}. 
\qed

\medskip

Let $I_{1,i}=B_1\cap f^i B$ for $i\ge 1$. 
Then $I_{1,1}=B_1\cap fB$ is a $\delta$-stable prime ideal of $B_1$. 
By Lemma \ref{Lemma 3.1}, we have a sequence of $G_a$-equivariant affine modifications 
$$A\subset B_0\subset B_1\subset B.$$
Note that 
$$\overline A \subset \overline B_0 \subset \overline B_1\subset B/fB$$
where $\overline B_1=B_1/I_{1,1}$. 
Let $\overline I_{1,i}$ be the image of $I_{1,i}$ by the surjection $B_1\to B_1/fB_1$.

\begin{lem}\label{Lemma 3.2}
The followings hold. 
\begin{enumerate}
\item[(1)] If $I_{1,1}=fB_1$, then $I_{1,i}=f^iB_1$ for all $i\ge 1$ and $B_1=B$. 
\item[(2)] Suppose that $I_{1,1}\supsetneq fB_1$. 
Let $g_1$ be an element of $I_{1,1}$ such that $\overline g_1 \in B_1/fB_1$ is a generator of 
$\overline I_{1,1}\otimes_{\overline B_0}Q(\overline B_0)=Q(\overline B_0)[\overline y_1]$
and write 
$$g_1=f^{\ell_1} y_2$$
where $\ell_1\ge 1$ and $y_2\in B\setminus fB$.
Let $B_2=B_1[f^{-\ell_1}I_{1,\ell_1}]$.
Then there exists $c_1 \in B_0 \setminus I_{0,1}$ such that
$$B_2[c_1^{-1}]=B_1[c_1^{-1}][y_2]$$
and 
$$(B_2/fB_2)[\overline c_1^{-1}]=\overline B_1[\overline c_1^{-1}][\overline y_2]
=(\overline B_1[\overline c_1^{-1}])^{[1]}.$$
\item[(3)] Suppose that $I_{1,1}\supsetneq fB_1$ and 
$\overline A$ and $\overline B_0$ are factorial. 
Take $g_1\in I_{1,1}$ so that $\overline g_1\in B_1/fB_1=\overline B_0[\overline y_1]$ 
is primitive over $\overline B_0$. Then $I_{1,1}=(f, g_1)B_1$ and 
$$B_2=B_1[f^{-\ell_1}I_{1,\ell_1}]=B_1[y_2]=A[z,y_1,y_2].$$ 
Hence $B_2$ is factorial and 
$$B_2/fB_2=\overline B_1[\overline y_2]=\overline B_1^{[1]}.$$
Further, $\overline{\partial_{y_1}g_1} \notin \overline I_{1,1}$, $p-\ell_0\ge \ell_1$ and 
$$\delta(y_2)=f^{p-\ell_0-\ell_1}g'\cdot \partial_{y_1}g_1+f^{p-\ell_0-\ell_1+1} b_1
\in B_1$$ 
for $b_1\in B_1$. Here, $\partial_{y_1}g_1=\frac{\partial g_1}{\partial y_1}(z,y_1)$. 
\end{enumerate}
\end{lem}
\Proof 
(1) We show $I_{1,i}=f^iB_1$ by induction on $i$. 
Let $h\in I_{1,i+1}$. Then $h$ is written as $h=f^{i+1}b$ for $b \in B$.
Since $h\in I_{1,i}=f^iB_1$, $h=f^ib_1$ for $b_1\in B_1$. 
Then it follows that $b_1=fb\in B_1\cap fB=I_{1,1}=fB_1$.
Hence $b\in B_1$ and $h\in f^{i+1}B_1$. Thus the first assertion follows. 

Since $B=A[z][f^{-\mu}I]$ for $I=A[z]\cap f^\mu B$, it follows that $B=B_1[f^{-\mu}I_{1,\mu}]$.
We have $B=B_1$ since $I_{1,\mu}=f^\mu B_1$. 

(2) Since $I_{1,1}\supsetneq fB_1$, $\overline I_{1,1}\neq (0)$.
Note that $\overline g_1$ has the minimal degree with respect to $\overline y_1$
in $\overline I_{1,1}$ and is irreducible in $B_1/fB_1$. 
Let $I_{1,1}=(f, h_1, h_2, \ldots,h_r)B_1$ for $h_1, \ldots, h_r \in B_1$. 
We may assume $h_1=g_1$ and $\{f, h_1=g_1, h_2, \ldots, h_r\}$ is a minimal set of generators.
We have $\overline a_i \overline h_i=\overline g_1\overline q_i$ in $B_1/fB_1$
for $a_i \in B_0\setminus I_{0,1}$ and $q_i \in B_1$. 
Let $c_1=c_0 a_2\cdots a_r \in B_0\setminus I_{0,1}$.
Then $c_1$ satisfies the required conditions. 

(3) The first and second assertions and the factoriality of $B_2$ follow
from the same argument as in Lemma \ref{Lemma 3.1} (5). 
Note that $\overline{\partial_{y_1}g_1}=\partial_{\overline y_1}\overline g_1$
is a nonzero element of $B_1/fB_1=\overline B_0[\overline y_1]$. 
Hence $\overline{\partial_{y_1}g_1}\notin \overline I_{1,1}
=\overline g_1 \overline B_0[\overline y_1]$.  

By $f^{\ell_1} y_2=g_1 \in B_0[y_1]=A[z,y_1]$, we have 
$$f^{\ell_1}\delta(y_2)=\partial_zg_1 \cdot \delta(z)+\partial_{y_1}g_1\cdot \delta(y_1)
=f^p\partial_zg_1+f^{p-\ell_0}g'(z)\partial_{y_1}g_1.$$
If $\ell_1>p-\ell_0$, then $\partial_{y_1}g_1\cdot g'(z)\in I_{1,1}$. 
Then $\overline{\partial_{y_1}g_1}\cdot \overline{g'}\in \overline I_{1,1}$,
which is a contradiction. 
Note that $\overline{g'}$ is nonzero as an element of
$\overline B_0 \subset B_1/fB_1=\overline B_0[\overline y_1]$. 
Thus $p-\ell_0\ge \ell_1$ and the assertion follows. 
\qed

\medskip

We inductively define $B_{i+1}$ for $i\ge 2$ by repeating the argument in Lemma \ref{Lemma 3.2}. 
Suppose that $B_j$ for $0 \le j \le i$ are defined by the same argument in Lemma \ref{Lemma 3.2}. 
Let $I_{i,j}=B_i\cap f^j B$ for $j\ge 1$. 
Then $I_{i,1}=B_i\cap fB$ is a $\delta$-stable prime ideal of $B_i$. 
Let $\overline I_{i,1}$ be the image of $I_{i,1}$ by the surjection $B_i\to B_i/fB_i$. 
For $0 \le j \le i$, $\overline B_j$ is defined by $\overline B_j=B_j/I_{j,1}$
where $I_{j,1}=B_j\cap fB$.   

\begin{lem}\label{Lemma 3.3}
The followings hold. 
\begin{enumerate}
\item[(1)] If $I_{i,1}=fB_i$, then $I_{i,j}=f^jB_i$ for all $j\ge 1$ and $B_i=B$. 
\item[(2)] Suppose that $I_{i,1}\supsetneq fB_i$. 
Let $g_i$ be an element of $I_{i,1}$ such that $\overline g_i \in B_i/fB_i$ is a generator of 
$\overline I_{i,1}\otimes_{\overline B_{i-1}}Q(\overline B_{i-1})
=Q(\overline B_{i-1})[\overline y_i]$ and write 
\begin{equation}
g_i=f^{\ell_i} y_{i+1} \label{(4)}
\end{equation}  
where $\ell_i\ge 1$ and $y_{i+1}\in B\setminus fB$.
Let $B_{i+1}=B_i[f^{-\ell_i}I_{i,\ell_i}]$.
Then there exists $c_i \in B_{i-1} \setminus I_{i-1,1}$ such that
$B_{i+1}[c_i^{-1}]=B_i[c_i^{-1}][y_{i+1}]$ and 
$$(B_{i+1}/fB_{i+1})[\overline c_i^{-1}]=\overline B_i[\overline c_i^{-1}][\overline y_{i+1}]
=(\overline B_i[\overline c_i^{-1}])^{[1]}.$$
\item[(3)] Suppose that $I_{i,1}\supsetneq fB_i$ and 
$\overline A$ and $\overline B_j$ for $0 \le j \le i-1$ are factorial.
Take $g_i \in I_{i,1}$ so that $\overline g_i\in B_i/fB_i=\overline B_{i-1}[\overline y_i]$ 
is primitive over $\overline B_{i-1}$. Then $I_{i,1}=(f,g_i)B_i$ and  
$$B_{i+1}=B_i[f^{-\ell_i}I_{i,\ell_i}]=B_i[y_{i+1}]=A[z,y_1,\ldots, y_{i+1}].$$ 
Hence $B_{i+1}$ is factorial and 
$$B_{i+1}/fB_{i+1}=\overline B_i[\overline y_{i+1}]=\overline B_i^{[1]}.$$ 
Further, $\overline{\partial_{y_i}g_i} \notin \overline I_{i,1}$,
$p-\ell_0 \cdots -\ell_{i-1}\ge \ell_i$ and 
\begin{align}
\delta(y_{i+1})&= \frac{\delta(g_i)}{f^{\ell_i}} \notag \\
&= f^{p-\ell_0 \cdots-\ell_i}g' \cdot \partial_{y_1}g_1 \cdots
\cdot \partial_{y_i}g_i+f^{p-\ell_0 \cdots-\ell_i+1} b_i \in B_i \label{(5)} 
\end{align}
for $b_i\in B_i$.
\end{enumerate}
\end{lem}
\Proof 
(1) and (3) follow from the same argument of the proof of Lemma \ref{Lemma 3.2}. 

(2) We have 
$$B_i[c_{i-1}^{-1}]=B_{i-1}[c_{i-1}^{-1}][y_i]$$
and 
$$(B_i/fB_i)[\overline c_{i-1}^{-1}]=\overline B_{i-1}[\overline c_{i-1}^{-1}][\overline y_i]
=(\overline B_{i-1}[\overline c_{i-1}^{-1}])^{[1]}$$
for $c_{i-1}\in B_{i-2}\setminus I_{i-2,1}$ and $y_i\in B_i\setminus fB$ such that
$f^{\ell_{i-1}}y_i=g_{i-1}\in I_{i-1,1}$. 
The assertion follows by the same argument as in Lemma \ref{Lemma 3.2} (2). 
\qed

\medskip

\begin{lem}\label{Lemma 3.4}
There exists a positive integer $\nu$ such that $B_\nu=B$.
\end{lem}
\Proof
Let $b_1, \ldots, b_r$ be the generators of $B$ over $\C$. 
Since $B[f^{-1}]=A[f^{-1}][z]$, 
there exists a nonnegative integer $\mu_i$ such that $f^{\mu_i}b_i\in A[z]$ for $1 \le i \le r$. 
If $\mu_i\le \ell_0$, then $f^{\ell_0}b_i\in I_{0,\ell_0}$ and $b_i\in B_1$.
If $\mu_i> \ell_0$, then $f^{\mu_i}b_i\in I_{0,\ell_0}$ and $f^{\mu_i-\ell_0}b_i\in B_1$.
If $\mu_i-\ell_0\le \ell_1$, then $f^{\ell_1}b_i\in I_{1,\ell_1}$ and $b_i\in B_2$. 
If $\mu_i-\ell_0 > \ell_1$, then $f^{\mu_i-\ell_0-\ell_1}b_i\in B_2$.
Continuing the argument, every $b_i$ is contained in $B_\nu$ for some $\nu >0$.
\qed

\medskip

\Remark We can take $\mu=\max \{\mu_i ; 1 \le i \le r\}$ in the equation (\ref{(2)}). 
In the case that $\beta \in A$ in the equation (\ref{(1)}) is not necessarily a unit, 
it holds that $B_\nu=B^{(1)}$ for some $\nu >0$. 

\bigskip

By Lemmas \ref{Lemma 3.3} and \ref{Lemma 3.4},
we have a sequence of $G_a$-equivariant affine modifications
\begin{equation}
A[z]=B_0\subset B_1\subset B_2 \subset \cdots \subset B_{\nu-1}\subset B_{\nu}=B \label{(6)}
\end{equation}
where $B_{i+1}[c_i^{-1}]=B_i[c_i^{-1}][y_{i+1}]$ for $c_i \in B_{i-1}\setminus I_{i-1,1}$
with the relation (\ref{(4)}), and
the associated sequence 
\begin{equation*}
X=X_\nu \xrightarrow{\sigma_\nu} X_{\nu-1} \xrightarrow{\sigma_{\nu-1}} \cdots
\to X_1 \xrightarrow{\sigma_1} X_0=Y\times \A^1 
\end{equation*}
of $G_a$-varieties and $G_a$-equivariant birational morphisms over $Y$
where $X_i=\Spec B_i$ for $0 \le i \le \nu$. 
If $\overline A$ and $\overline B_i$ for $0 \le i \le \nu-2$ are factorial,
then $B_{i+1}=B_i[y_{i+1}]$ for $0 \le i < \nu$ and the $G_a$-action is given by (\ref{(5)}).
Hence in this case $X=\Spec B$ is isomorphic to a factorial affine variety over $Y$ defined 
in $Y\times \A^{\nu+1}$ by the equations 
\begin{align}
f^{\ell_0} y_1 & =g(z), \notag \\ 
f^{\ell_1} y_2 & =g_1(z,y_1), \notag \\ 
\cdots & \cdots \cdots \qquad \quad , \\ \label{(7)}
f^{\ell_{\nu-2}} y_{\nu-1} & =g_{\nu-2}(z,y_1, \ldots, y_{\nu-2}), \notag \\ 
f^{\ell_{\nu-1}} y_\nu & =g_{\nu-1}(z, y_1, \ldots, y_{\nu-1})  \notag
\end{align}    
where $p-\ell_0 \cdots -\ell_{i-1}\ge \ell_i>0$ for $0 \le i < \nu$.    
Since $\delta$ is irreducible, it follows that $p=\ell_0 +\cdots +\ell_{\nu-1}$ by (\ref{(5)}). 
Note that
$$\overline A \subset \overline B_0 \subset \cdots \subset \overline B_{\nu-1}\subset
\overline B_{\nu}=B/fB$$
where $\overline B_i=B_i/I_{i,1}$ and $I_{i,1}=B_i\cap fB$ for $0 \le i \le \nu$.  

Let $\overline \delta$ be the lnd on $B/fB$ induced by $\delta$ and
let $Z=\Spec \Ker \overline \delta$. 
Note that $\overline \delta\neq 0$ since $\delta$ is irreducible. 
Let $\overline \pi : E=\Spec B/fB \to Z$ be the quotient morphism. 
Then the restriction $\pi|_E : E \to V_Y(f)=\Spec \overline A$ splits as
$\pi|_E=\tau\circ \overline \pi$ where $\tau : Z \to V_Y(f)$ is the morphism
induced by the inclusion $\overline A \hookrightarrow \Ker \overline \delta$. 
Note that $\Ker \overline \delta$ is algebraic over $\overline A$
since the transcendental degree of $Q(\Ker \overline \delta)$ and that of $Q(\overline A)$
are equal. 

We have the following. 

\begin{thm}\label{Theorem 3.5}
Let $X=\Spec B$ be a factorial affine variety with a $G_a$-action 
corresponding to an irreducible lnd $\delta$. 
Suppose that $A=\Ker \delta$ is an affine $\C$-domain and 
the quotient morphism $\pi: X \to Y=\Spec A$ is faithfully flat. 
Then the following assertions are equivalent.
\begin{enumerate}
\item[(1)] The quotient morphism $\pi$ is a trivial $\A^1$-bundle. 
\item[(2)] There exists a closed subset $V(\alpha)$ of $Y$ of codimension one satisfying  
\begin{enumerate}
\item[(i)] every closed fiber of $\pi$ over $Y\setminus V(\alpha)$ is $\A^1$, 
\item[(ii)] for every prime factor $f_i$ of $\alpha$ in $A$, 
the morphism $\tau_i : Z_i=\Spec \Ker \delta_i \to V_Y(f_i)=\Spec A/f_iA$ induced by
the inclusion $A/f_iA \hookrightarrow \Ker \delta_i$ is birational
where $\delta_i$ is the lnd on $B/f_i B$ induced by $\delta$.  
\end{enumerate}  
\end{enumerate}
\end{thm}
\Proof
If (1) holds, there exists a slice $s\in B$ of $\delta$ and $B=A[s]$. 
Then (2) holds obviously. 

We show (2) implies (1). 
Since $B$ is faithfully flat over $A$, the plinth ideal $\pl(\delta)$ is principal.
If $\pl(\delta)$ is a unit ideal, then $B=A^{[1]}$ and the assertion follows. 
Suppose $\pl(\delta)\neq A$. 
By the condition (i), we have $\pl(\delta)=f_1^{p_1}\cdots f_s^{p_s} A$
where $p_i\ge 0$ and $f_1, \ldots, f_s$ are distinct prime factors of $\alpha$. 
We may assume $p_i>0$ for every $i$. 
Let $z\in B$ be a local slice such that $\delta(z)=f_1^{p_1}\cdots f_s^{p_s}$. 
Take any $f_i$, say $f_1$, and let $f=f_1$ and set $\beta=f_2^{p_2}\cdots f_s^{p_s}$. 
The lnd $\delta$ extends to an lnd $\delta_\beta$ on $B_\beta=B[\beta^{-1}]$ and
the plinth ideal of $\delta_\beta$ is generated by $f^{p_1}$.
Let $\overline \delta$ be the lnd on $B_\beta/fB_\beta$ induced by $\delta_\beta$.  
We have 
$$(A/fA)_{\overline \beta}=A_\beta/fA_\beta \subset A_\beta[z]/(A_\beta[z]\cap fB_\beta) 
\subset \Ker \overline \delta \subset B_\beta/fB_\beta.$$ 
By Lemma \ref{Lemma 3.1} (3), we have
$$[Q(\Ker \overline \delta): Q(A/fA)] \ge [Q(\overline B_0) : Q(A/fA)]>1$$
where $\overline B_0=A[z]/(A[z]\cap fB)$. 
This is a contradiction to the condition (ii). 
Hence the assertion follows. 
\qed

\medskip

\Remarks (1) The restriction $\pi|_{V(f_i)} : V(f_i)=\Spec B/f_i B \to V_Y(f_i)$ splits as
$\pi|_{V(f_i)}=\tau_i\circ \pi_i$ where $\pi_i: V(f_i) \to Z_i$ is the quotient morphism.
Since $\delta$ is irreducible, $\delta_i\neq 0$ and a general fiber of $\pi_i$ is $\A^1$. 
Hence a general fiber of $\pi|_{V(f_i)}$ consists of a disjoint union
of $m_i$ affine lines where $m_i=[Q(\Ker \delta_i): Q(A/f_iA)]>1$ by Theorem \ref{Theorem 2.2}.

(2) The conditions (i) and (ii) of Theorem \ref{Theorem 3.5} (2) imply 
$\codim_Y\Sing (\pi)>1$ (cf. the remark of Proposition \ref{Proposition 2.5}). 

\bigskip

In the case that the $G_a$-action is free, Theorem \ref{Theorem 3.5} implies the following. 

\begin{thm}\label{Theorem 3.6}
Let $X=\Spec B$ be a smooth factorial affine variety with a free $G_a$-action 
corresponding to an lnd $\delta$. 
Suppose that $A=\Ker \delta$ is an affine $\C$-domain and 
the quotient morphism $\pi : X \to Y=\Spec A$ is equi-dimensional and surjective. 
Then the following assertions are equivalent. 
\begin{enumerate}
\item[(1)] The quotient morphism $\pi$ is a trivial $\A^1$-bundle. 
\item[(2)] There exists a closed subset $V(\alpha)$ of $Y$ of codimension one satisfying  
\begin{enumerate}
\item[(i)] every closed fiber of $\pi$ over $Y\setminus V(\alpha)$ is $\A^1$, 
\item[(ii)] for every prime factor $f_i$ of $\alpha$ in $A$,
the morphism $\tau_i : Z_i=\Spec \Ker \delta_i \to V_Y(f_i)=\Spec A/f_iA$ induced by
the inclusion $A/f_iA \hookrightarrow \Ker \delta_i$ is birational
where $\delta_i$ is the lnd on $B/f_i B$ induced by $\delta$. 
\end{enumerate}  
\end{enumerate}  
\end{thm}
\Proof
By Proposition \ref{Proposition 2.5}, $\pi$ is faithfully flat. 
Hence the assertion follows from Theorem \ref{Theorem 3.5}. 
\qed

\medskip

As observed below, it rarely happens that the quotient morphism is a trivial $\A^1$-bundle.

Let $A$ be anew an affine $\C$-domain and
let $I$ be the ideal of $A[Z, Y_1,\ldots, Y_\nu]=A^{[\nu+1]}$.

First, consider the case that $I$ is generated by 
\begin{align*}
& x^{\ell_0} Y_1 -g(Z),  \\ 
& x^{\ell_1} Y_2 -g_1(Z,Y_1),  \\ 
& \cdots  \cdots \cdots \qquad \quad , \\
& x^{\ell_{\nu-2}} Y_{\nu-1} -g_{\nu-2}(Z,Y_1, \ldots, Y_{\nu-2}), \\ 
& x^{\ell_{\nu-1}} Y_\nu -g _{\nu-1}(Z, Y_1, \ldots, Y_{\nu-1})  
\end{align*}    
where $x\in A$ is a prime element, $g(Z)\in A[Z]\setminus A$, $\ell_i>0$ for $0 \le i < \nu$ and
$g_i(Z,Y_1,\ldots,Y_i)\in A[Z,Y_1,\ldots, Y_i]\setminus A[Z,Y_1,\ldots, Y_{i-1}]$
for $0 < i < \nu$. 
Let
$$B=A[Z, Y_1,\ldots, Y_\nu]/I=A[z, y_1,\ldots, y_\nu]$$
where $z$ and $y_i$ denote the residue class of $Z$ and $Y_i$, respectively. 
Then $B[x^{-1}]=A[x^{-1}][z]=A[x^{-1}]^{[1]}$. 
Let $\pi : X=\Spec B \to Y=\Spec A$ be the morphism defined by the inclusion $A\hookrightarrow B$.
Then $\pi$ is an $\A^1$-fibration. 

\begin{prop}\label{Proposition 3.7}
Let $B$ be as above. 
Suppose that $B$ is factorial and faithfully flat over $A$ and $A$ is factorially closed in $B$. 
If $z-a\notin xB$ for any $a\in A$ and $\partial_{y_{\nu-1}}g_{\nu-1}\notin x B$, 
then a general fiber of $\pi|_{V(x)} : V(x)=\Spec B/xB \to V_Y(x)=\Spec A/xA$ consists of
a disjoint union of $m$ affine lines for $m>1$.  
Hence $\pi$ is not a trivial $\A^1$-bundle. 
\end{prop}
\Proof
Note that $A$ is factorial and $x$ is prime in $B$ since $A$ is factorially closed in $B$. 
Note also that $g(Z)\notin xA[Z]$ and $g_i(Z,Y_1,\ldots,Y_i)\notin xA[Z,Y_1,\ldots, Y_i]$
for $0 < i < \nu$ since $B$ is an integral domain. 
Let $p=\ell_0+\cdots +\ell_{\nu-1}$. 
We define a locally nilpotent $A$-derivation $\delta$ on $B$ by
$$\delta(z)=x^p, \quad \delta(y_1)=x^{p-\ell_0} g'(z), \quad 
\delta(y_i)=\frac{\delta(g_{i-1})}{x^{\ell_{i-1}}}$$
for $2\le i \le \nu$. Note that $\delta(y_i)\in x^{p-\ell_0 \cdots -\ell_{i-1}} B$.
Since $\partial_{y_{\nu-1}}g_{\nu-1}\notin x B$, $\delta$ is irreducible. 
We have $\Ker \delta=A$. In fact, since $B[x^{-1}]=A[x^{-1}][z]$ and
the kernel of the extended lnd of $\delta$ onto $B[x^{-1}]$ is $A[x^{-1}]$,
$\Ker \delta=B\cap A[x^{-1}]$. 
By the factorially-closedness of $A$ in $B$, it follows that $B\cap A[x^{-1}]=A$ and 
we have $\Ker \delta=A$.
Since $B$ is faithfully flat over $A$, the plinth ideal $\pl(\delta)$ is principal and
$\pl (\delta)=x^e A$ for $0 \le e \le p$.
Suppose $e=0$, i.e., $\delta$ has a slice $s\in B$. 
Then $\delta(z-x^p s)=0$, and $z-a=x^p s\in xB$ for some $a\in A$.
This contradicts to the assumption $z-a\notin xB$ for any $a\in A$. 
Hence $\pl (\delta)=x^e A$ for $e>0$. 
Let $\overline \delta$ be the lnd on $B/xB$ induced by $\delta$ and
let $\overline \pi : V(x)\to \Spec \Ker \overline \delta$ be the quotient morphism.
Let $\tau : \Spec \Ker \overline \delta \to V_Y(x)=\Spec A/xA$ be the morphism induced by
$A/xA \hookrightarrow \Ker \overline \delta$. 
We have $\pi|_{V(x)}=\tau\circ \overline \pi$.
Using Lemma \ref{Lemma 3.1} (3) and
the same argument as in the proof of Theorem \ref{Theorem 3.5}, we have 
$m=[Q(\Ker \overline \delta): Q(A/xA)]>1$. Hence $\tau$ is not birational 
and $\pi$ is not a trivial $\A^1$-bundle by Theorem \ref{Theorem 3.5}. 
Since a general fiber of $\overline \pi$ is $\A^1$, 
a general fiber of $\pi|_{V(x)}$ is a disjoint union of $m$ affine lines
as observed in the remark of Theorem \ref{Theorem 3.5}. 
\qed 

\medskip

\Remark Let $x_1, \ldots, x_s\in A$ be distinct prime elements 
and $g(Z)$ and $g_i(Z,Y_1,\ldots,Y_i)$ be as above. 
Proposition \ref{Proposition 3.7} holds true for $B=A[Z, Y_1,\ldots, \\
Y_\nu]/I$ 
where $I$ is generated by $x_1^{\ell_{0,1}}\cdots x_s^{\ell_{0,s}} Y_1-g(Z)$ and 
$x_1^{\ell_{i,1}}\cdots x_s^{\ell_{i,s}} Y_{i+1} -g_i(Z,Y_1, \ldots, Y_i)$ for $1 \le i < \nu$
where $\ell_{i,j}>0$ for $0 \le i < \nu$ and $1 \le j \le s$.
Precisely, for such $B$ which is factorial and faithfully flat over $A$ and
in which $A$ is factorially closed, 
if for every $1 \le i \le s$, $\partial_{y_{\nu-1}}g_{\nu-1}\notin x_i B$ and
$z-a\notin x_iB$ for any $a\in A$, 
then a general fiber of $\pi|_{V(x_i)}$ is a disjoint union of $m_i$ affine lines with
$m_i=[Q(\Ker \delta_i): Q(A/x_iA)]>1$ where $\delta_i$ is the induced lnd on $B/x_iB$, 
and $\pi$ is not a trivial $\A^1$-bundle. 

\bigskip

Next, consider the case $\nu=1$.

Let $A$ be factorial and let $\alpha=x_1^{p_1}\cdots x_s^{p_s}$
where $x_1,\ldots, x_s$ are distinct prime elements of $A$ and $p_i>0$ for $1 \le i \le s$. 
Let $B=A[Y,Z]/(\alpha Y-g(Z))$ where $g(Z)\in A[Z]\setminus A$ is an irreducible polynomial
such that $(x_i, g(Z))A[Z]$ is a prime ideal and $Z-a\notin (x_i, g(Z))A[Z]$ for every $i$.
Let $\pi : X=\Spec B \to Y=\Spec A$ be the $\A^1$-fibration 
defined by the inclusion $A\hookrightarrow B$.

\begin{prop}\label{Proposition 3.8}
Let $A$ and $B$ be as above. Suppose that $B$ is faithfully flat over $A$.  
If $g'(Z)\notin (x_i, g(Z))A[Z]$ for $1 \le i \le s$, 
then a general fiber of $\pi|_{V(x_i)}$ is a disjoint union of $m_i$ affine lines with $m_i>1$ 
for every $i$ and $\pi$ is not a trivial $\A^1$-bundle. 
\end{prop}
\Proof
Note that each $x_i$ is prime in $B$ 
since $B/x_i B=\left(A[Z]/(x_i, g(Z))A[Z]\right)^{[1]}$ is an integral domain. 
Hence $B$ is factorial since $B[\alpha^{-1}]=A[\alpha^{-1}][z]$ \cite{Nag}. 
We define a locally nilpotent $A$-derivation $\delta$ on $B$ by
$\delta(z)=\alpha$ and $\delta(y)=g'(z)$. 
Then $\delta$ is irreducible by the assumption $g'(Z)\notin (x_i, g(Z))A[Z]$ for $1 \le i \le s$.
The element $z\in B$ is a local slice of $\delta$ and $\Ker \delta=A$. 
Since $B$ is faithfully flat over $A$, the plinth ideal $\pl(\delta)$ is principal and
$\pl (\delta)=x_1^{e_1}\cdots x_s^{e_s} A$ for $0 \le e_i \le p_i$. 
By a similar argument in the proof of Proposition \ref{Proposition 3.7}, $e_i>0$ for every $i$.
In fact, suppose $e_i=0$ for some $i$. 
There exists $b\in B$ such that $\delta(b)=x_1^{e_1}\cdots x_s^{e_s}$. 
Since $\delta(z-x_1^{p_1-e_1}\cdots x_i^{p_i}\cdots x_s^{p_s-e_s}b)=0$,
$z-a\in x_i B$ for some $a\in A$.
This is a contradiction to the assumption $Z-a\notin (x_i, g(Z))A[Z]$ for $1 \le i \le s$. 
Hence $e_i>0$ for every $i$.
Let $\delta_i$ be the induced lnd on $B/x_i B$.  
We have $\pi|_{V(x_i)}=\tau_i\circ \pi_i$ where $\pi_i : V(x_i)\to \Spec \Ker \delta_i$ is the
quotient morphism and $\tau_i : \Spec \Ker \delta_i \to \Spec A/x_i A$ is the morphism
defined by $A/x_i A\hookrightarrow \Ker \delta_i$. 
By the same argument as in the proof of Proposition \ref{Proposition 3.7}, 
a general fiber of $\pi|_{V(x_i)}$ is a disjoint union of $m_i$ affine lines with
$m_i=[Q(\Ker \delta_i): Q(A/x_iA)]>1$, and $\pi$ is not a trivial $\A^1$-bundle.  
\qed

\medskip

\Remarks (1) It holds that $e_i=p_i$ for every $i$. In fact, if $0<e_i<p_i$ for some $i$.
Then $\delta(z-x_1^{p_1-e_1}\cdots x_i^{p_i-e_i}\cdots x_s^{p_s-e_s}b)=0$ implies 
$z-a\in x_i B$ for some $a\in A$. This is a contradiction to the assumption.
Hence $e_i=p_i$ for every $i$ and $\pl(\delta)=\alpha A$. 

(2) Let $A$ be a factorial affine domain and $B=A[Y,Z]/(\alpha Y-g(Z))$
where $A[Y,Z]=A^{[2]}$, $\alpha \in A$ and $g(Z)\in A[Z]\setminus A$ is an irreducible polynomial.
Then the factoriality of $B$ is studied in \cite{DFN}.

If $B\cong A^{[1]}$, $\alpha Y-g(Z)$ is a coordinate of $B$ over (any $\Q$-algebra) $A$ by
\cite{ER}. If $A=k[X_1,\ldots, X_m,T]=k^{[m+1]}$ for a field $k$ of characteristic $0$
and $g(Z)=f(Z,T)+h(X_1, \ldots,X_m,T,Z)$ with $f(Z,T)\in k[Z,T]$ and
$h(X_1, \ldots,X_m,T,Z)\in k[X_1,\ldots, X_m,T]$ which is divided by
every prime factor of $\alpha$, criteria for $B=k^{[m+2]}$ are given under some assumptions by
\cite{GhoGupPal}.

\bigskip
\medskip

\section{Smooth acyclic affine varieties with free $G_a$-actions}

First, we consider the case that $X$ is a smooth acyclic affine threefold
with an $\A^2$-fibration $f : X \to \A^1$ and 
give a short proof to a characterization of $\A^3$ due to
Miyanishi \cite{Miy4}, Kaliman \cite{Kal2}.

\begin{thm}\label{Theorem 4.1}
Let $X=\Spec B$ be a smooth acyclic affine threefold. 
If there exists an $\A^2$-fibration $f : X \to \A^1$ such that 
every closed fiber of $f$ is factorial, then $X\cong \A^3$ and $f$ is a variable of $B=\C^{[3]}$. 
\end{thm}
\Proof
By Fujita \cite{Fuj}, $X$ is factorial and $B^*=\C^*$. 
The existence of an $\A^2$-fibration $f : X \to \A^1$ implies that 
there exists $\beta\in \C[f]\subset B$ such that $B[\beta^{-1}]=\C[f, \beta^{-1}, x,y]$,
a polynomial ring in two variables over $\C[f,\beta^{-1}]$.
Let $D$ be a $\C[f,\beta^{-1}]$-derivation on $B[\beta^{-1}]$ defined by $D(x)=1$ and $D(y)=0$.
Then the derivation $\delta=\beta^m D$ restricts to an lnd on $B$ for some integer $m \ge 0$
since $B$ is finitely generated over $\C$. 
We may assume $\delta$ is irreducible. 
Let $A=\Ker \delta$. Then $A$ is factorial and $A^*=B^*=\C^*$.
Further, $f\in A$ and $\delta(x)$ is a product of some factors of $\beta$. 
Hence the $\A^2$-fibration $f : X \to \A^1$ is equivariant
with respect to the $G_a$-action associated with $\delta$ and $A[\beta^{-1}]=\C[f,\beta^{-1},y]$. 
The $G_a$-equivariant $\A^2$-fibration $f : X \to \A^1$ 
splits as a composite of the quotient morphism $\pi : X \to Y$ and $f_Y: Y \to \A^1$.  
Since $f_Y : Y \to \A^1$ is an $\A^1$-fibration, $A=\C^{[2]}$ 
by a characterization of affine plane \cite{Miy} and $f$ is a variable of $A$.  
For every $c\in \C$, the fiber $W_c=\Spec B/(f-c)B$ of $f$ is factorial and
has a non-trivial $G_a$-action induced by $\delta$.
The lnd $\delta_c$ on $B/(f-c)B$ induced by $\delta$ is written as
$\delta_c=a_c d_c$ for an irreducible lnd $d_c$ on $B/(f-c)B$ and $a_c\in \Ker d_c$.
Consider the $G_a$-action on $W_c$ associated with $d_c$.
By Lemma \ref{Lemma 2.1}, $W_c$ is smooth. 
Further, by Proposition \ref{Proposition 2.6}, 
each $W_c$ is acyclic, hence $W_c\cong \A^2$ by Miyanishi \cite{Miy3}.
By \cite[Theorem 2.1.16]{GMM3}, \cite{Sat},
$f: X\to \A^1$ is an $\A^2$-bundle, and hence $X\cong \A^3$ and $f$ is a variable. 
\qed

\medskip

\Remark
The existence of an $\A^2$-fibration $f : X \to \A^1$ is equivalent to
the existence of an $\A^2$-cylinder $U\cong C\times \A^2$ in $X$ for an affine curve $C$
(\cite{Miy4}, \cite{Kal2}).
The above proof of Theorem \ref{Theorem 4.1} shows that
if a smooth acyclic affine threefold $X$ contains an $\A^2$-cylinder $U$
and each irreducible component of $X\setminus U$ is factorial, then $X\cong \A^3$. 

\medskip

Let $X$ be as in Theorem \ref{Theorem 4.1} and
suppose that $X$ has a $G_a$-action and a $G_a$-equivariant $\A^2$-fibration $f : X \to \A^1$ 
where $G_a$ acts trivially on $\A^1$. 
Theorem \ref{Theorem 4.1} does not imply that
the quotient morphism $\pi :X \to Y=\A^2$ is an $\A^1$-bundle
as the following simple example shows.

\medskip

\noindent{\bf Example 4.1}\label{Example 4.1} \: 
Let $B=\C[x,y,z]$ be a polynomial ring and let $X=\Spec B$.
Define an lnd $\delta$ on $B$ by $\delta(x)=0$, $\delta(y)=-2z$, $\delta(z)=x^2$.
Then $A=\Ker \delta=\C[x,t]=\C^{[2]}$ where $t=x^2y+z^2$ and let $Y=\Spec A$.
The quotient morphism $\pi : X \to Y$ is faithfully flat with $\Sing (\pi)\neq \emptyset$. 
In fact, $\Sing (\pi)$ is the closed set $V(x)\cong \A^1$ of $Y=\A^2_{(x,t)}$. 
For $Q=(0, \beta)\in V(x)$, the fiber $\pi^*(Q)$ is $\A^1\amalg \A^1$ if $\beta\neq 0$
and $\A^1$ with multiplicity $2$ if $\beta=0$. 
The fixed point locus consists of the fiber $\pi^{-1}(O)$ for $O=(0,0)\in Y$.
The projection $x : X \to \A^1$ is an trivial $\A^2$-bundle and $x$ is a variable of $B$, 
however, $\pi : X \to Y$ is not an $\A^1$-bundle. 

\medskip

Let $X=\Spec B$ be an affine variety of dimension $\ge 3$ with a $G_a$-action and 
a $G_a$-equivariant fibration $f : X \to \A^1$ where $G_a$ acts trivially on $\A^1$. 
Let $\pi : X \to Y=\Spec A$ be the quotient morphism with $A=\Ker \delta$
where $\delta$ is the lnd on $B$ associated with the $G_a$-action.
Since $f : X \to \A^1$ is $G_a$-equivariant, $f\in A$ defines the morphism $f_Y : Y \to \A^1$ 
and $f$ splits as $X\stackrel{\pi}{\to} Y \stackrel{f_Y}{\to} \A^1$.
For $c\in \C$, let $Y_c=f_Y^{-1}(c)$ and $W_c=f^{-1}(c)$.
The $G_a$-action restricts to $W_c=\Spec B/(f-c)B$.
Let $\pi_c: W_c\to Z_c=\Spec \Ker \delta_c$ be the quotient morphism
where $\delta_c$ is the induced lnd on $B/(f-c)B$. 
Then the restriction $\pi|_{W_c} : W_c \to Y_c$ splits as $\pi|_{W_c}=\nu_c\circ \pi_c$
where $\nu_c : Z_c \to Y_c=\Spec A/(f-c)A$ is the morphism defined by the inclusion
$A/(f-c)A \hookrightarrow \Ker \delta_c$. 
As for the singular locus $\Sing (\pi)$, we have the following result.  

\begin{lem}\label{Lemma 4.2}
Let $X=\Spec B$ be a smooth factorial affine variety of dimension $n\ge 3$
with a free $G_a$-action. Let $\pi : X \to Y$ be the quotient morphism. 
We assume that $Y=\Spec A$ is an affine variety and $\pi$ is surjective and equi-dimensional. 
Suppose that there exists a $G_a$-equivariant fibration $f : X \to \A^1$
with a trivial $G_a$-action on $\A^1$ such that
a general closed fiber $F$ of $f$ is $G_a$-equivariantly isomorphic to
$Z \times \A^1$ where $G_a$ acts on an affine variety $Z$ trivially and 
on the second factor $\A^1$ by translation.
Then $\Sing (\pi)=\emptyset$ or
$\overline{\Sing (\pi)}$ is a union of finite fibers of $f_Y : Y\to \A^1$. 
\end{lem}
\Proof
Let $\delta$ be the lnd on $B$ associated with the $G_a$-action. 
Then $\delta$ is irreducible since the action is free, and
$A=\Ker \delta$ is an affine $\C$-domain by the assumption. 
By Proposition \ref{Proposition 2.5}, $Y=\Spec A$ is a smooth factorial affine variety and 
$\Sing(\pi)$ is of pure codimension one unless $\Sing (\pi)=\emptyset$. 
Suppose that $\Sing (\pi)\neq \emptyset$. 
As remarked above, $f\in A$ defines the morphism $f_Y : Y \to \A^1$ and $f=f_Y\circ \pi$. 
Choose $\lambda\in \C$ so that $Y_\lambda=\Spec A/(f-\lambda)A$ and
$W_\lambda=\Spec B/(f-\lambda)B$ are general fibers of $f_Y: Y \to \A^1$ and $f : X \to \A^1$,
respectively. 
Then $W_\lambda=\pi^{-1}(Y_\lambda)$ and $Y_\lambda$ are smooth by Bertini's theorem. 
Let $\delta_\lambda$ be the lnd on $B/(f-\lambda)B$ induced by $\delta$ and 
let $Z_\lambda=\Spec \Ker \delta_\lambda$. 
By the assumption, $\delta_\lambda$ has a slice, hence $Z_\lambda$ is an affine variety. 
The restriction $\pi|_{W_\lambda}: W_\lambda \to Y_\lambda$ 
is a composite of the quotient morphism $\pi_\lambda : W_\lambda \to Z_\lambda$ 
and a surjective morphism $\nu_\lambda : Z_\lambda \to Y_\lambda$ which is defined by
the inclusion $A/(f-\lambda)A \hookrightarrow \Ker \delta_\lambda$. 
Note that $Q(\Ker \delta_\lambda)$ is an algebraic extension of $Q(A/(f-\lambda)A)$.  
Since the $G_a$-action is free,
every closed fiber of $\pi$ is a disjoint union of finite $G_a$-orbits. 
Further, since $W_\lambda$ is $G_a$-equivariantly isomorphic to $Z_\lambda\times \A^1$ and
$\pi_\lambda : W_\lambda \to Z_\lambda$ is a trivial $\A^1$-bundle,  
every fiber of $\nu_\lambda$ consists of finite points. 
Since $Y_\lambda$ is a general fiber of $f_Y : Y \to \A^1$, 
$Y_\lambda \not\subset \overline{\Sing(\pi)}$, which implies that
a general fiber of $\pi |_{W_\lambda}$ is $\A^1$ (cf. \cite{Ma2}). 
Hence $\nu_\lambda$ is birational, 
and $Z_\lambda$ coincides with $Y_\lambda$ by Zariski's Main Theorem. 
Then $\pi|_{W_\lambda}$ coincides with the quotient morphism 
$\pi_\lambda : W_\lambda\to Z_\lambda$ which is a trivial $\A^1$-bundle by the assumption. 
Thus $Y_\lambda\cap \Sing (\pi)=\emptyset$ for general $\lambda$ and $\overline{\Sing (\pi)}$ is 
a union of finite fibers $Y_{c_i}=\Spec A/(f-c_i)A$ for $c_i\in \C$ and $1 \le i \le s$. 
\qed

\medskip

By Theorem \ref{Theorem 4.1} and Lemma \ref{Lemma 4.2}, we obtain the following result.

\begin{thm}\label{Theorem 4.3}
Let $X=\Spec B$ be a smooth acyclic affine fourfold with a free $G_a$-action.
Let $\pi :X \to Y$ be the quotient morphism where we assume $Y=\Spec A$ is an affine variety.   
If 
\begin{enumerate}
\item[(1)] $\pi$ is equi-dimensional and surjective, 
\item[(2)] $Y$ is acyclic,
\item[(3)] there exists a $G_a$-equivariant $\A^3$-fibration $f : X \to \A^1$ 
where the $G_a$-action on $\A^1$ is trivial such that 
every closed fiber of the induced morphism $f_Y : Y\to \A^1$ is factorial, 
\end{enumerate}
then $Y\cong \A^3$ and $f$ is a variable of $A=\C^{[3]}$.  

Further, suppose that
\begin{enumerate}
\item[(4)] the morphism $\nu_c : Z_c=\Spec \Ker \delta_c \to Y_c=\Spec A/(f-c)A$ defined by the 
inclusion $A/(f-c)A \hookrightarrow \Ker \delta_c$ is birational for every $c\in \C$.
Here, $\delta_c$ is the induced lnd on $B/(f-c)B$.   
\end{enumerate}  
Then $\pi : X \to Y$ is a trivial $\A^1$-bundle and $X\cong \A^4$.   
\end{thm}
\Proof
Since $X=\Spec B$ is smooth and acyclic, $B$ is factorial and $B^*=\C^*$. 
Let $\delta$ be an irreducible lnd on $B$ associated with the $G_a$-action and 
let $A=\Ker \delta$ which is an affine $\C$-domain by the assumption.
By Proposition \ref{Proposition 2.5}, $Y=\Spec A$ is a smooth factorial affine threefold. 
The quotient morphism $\pi : X \to Y$ is faithfully flat by the condition (1). 
The morphism $f_Y : Y \to \A^1$ induced by $f$ is a surjection. 
In fact, if $f^{-1}(c)=\emptyset$ for some $c\in \C$, then the ideal $(f-c)A$ is a unit ideal,
which implies $f-c\in A^*=B^*=\C^*$. This is a contradiction. 
Hence $f_Y: Y \to \A^1$ is surjective and so is $f$. 
For $c\in \C$, let $Y_c=\Spec A/(f-c)A$ and $W_c=\pi^{-1}(Y_c)=\Spec B/(f-c)B$. 
As observed above, $\pi |_{W_c}=\nu_c \circ \pi_c$
where $\pi_c : W_c \to Z_c$ is the quotient morphism. 
For a general $c\in \C$, $W_c\cong \A^3$ and $G_a$ acts freely on it. 
Hence for a general $c\in \C$, 
$Z_c\cong \A^2$ and $W_c\cong Z_c\times \A^1$ by Kaliman \cite{Kal1}. 
By the same argument as in the proof of Lemma \ref{Lemma 4.2}, 
$Y_c=Z_c\cong \A^2$ for a general $c\in \C$ and $f_Y : Y \to \A^1$ is an $\A^2$-fibration.
Since $Y_c$ is factorial for every $c \in \C$ and $Y$ is acyclic by the conditions (2) and (3), 
$Y\cong \A^3$ and $f$ is a variable of $A=\C^{[3]}$ by Theorem \ref{Theorem 4.1}. 

By Lemma \ref{Lemma 4.2}, we have $\Sing (\pi)\neq \emptyset$ or
$\overline{\Sing(\pi)}$ consists of finite fibers of $f_Y$.
If the condition (4) is satisfied,
then $\pi$ is a trivial $\A^1$-bundle by Theorem \ref{Theorem 3.6} and the assertion follows.
\qed

\medskip

If $G_a$ acts non-trivially on $X=\A^4$ and preserves a coordinate $x$, i.e., $\delta(x)=0$
for the lnd $\delta$ associated with the $G_a$-action on $X$, 
$A=\Ker \delta$ is an affine $\C$-domain by Bhatwadekar and Daigle \cite{BhD} 
and $Y=\Spec A$ is an affine threefold. 
By considering the projection $x : X \to \A^1$ as the $\A^3$-fibration $f$ 
in Theorem \ref{Theorem 4.3}, we obtain the following result of a free $G_a$-action on $\A^4$ 
preserving a coordinate.  

\begin{cor}\label{Corollary 4.4}
Let $G_a$ act freely on $X=\A^4$ preserving a coordinate $x$ of $X$. 
Suppose that the quotient morphism $\pi : X \to Y$ is equi-dimensional and surjective and 
$Y$ is acyclic. If every closed fiber of the induced morphism $x : Y \to \A^1$ is factorial,
then $Y\cong \A^3$ and $x$ is a coordinate of $Y$. 
\end{cor}

\medskip

\noindent{\bf Example 4.2}\label{Example 4.2} \: 
Let $X_0$ be a Koras-Russell $3$-fold defined by $X^2Y-g(X,Z,T)=0$ in $\A^4_{(X,Y,Z,T)}$
where $g(X,Z,T)=X+Z^2+T^3$ and let $\tilde X=X_0\times \A^1$.
Then $\tilde X$ is a smooth acyclic affine fourfold. 
Let $B=\C[X,Y,Z,T,W]/(X^2Y-g(X,T,Z))$ where $\C[X,Y,Z,T,W]$ is a polynomial ring in $5$ variables
and we regard $B$ as the coordinate ring of $\tilde X$. 
We denote the residue class of $X,Y,Z,T,W$ by $x,y,z,t,w\in B$, respectively. 
Define an lnd $\delta$ on $B$ by 
$$\delta(x)=\delta(t)=0, \quad \delta(y)=2z, \quad \delta(z)=x^2, \quad \delta(w)=1.$$ 
Since $\delta$ has a slice $w$, the $G_a$-action associated with $\delta$ on $\tilde X$ is free
and $B=A[w]$ where $A=\Ker \delta$. 
We have
$$A=\C[x,u,v,t]\cong \C[X,U,V,T]/(X^2U-g(X,V,T))$$
where $u=y-2zw+x^2w^2$ and $v=z-x^2w$.
Hence $Y=\Spec A \cong X_0$ is a smooth acyclic threefold. 
The function $x\in B$ defines a $G_a$-equivariant $\A^3$-fibration $x : \tilde X \to \A^1$.
The induced morphism $x_Y : Y \to \A^1$ is an $\A^2$-fibration with a single singular fiber 
$x_Y^{-1}(0)\cong C\times \A^1$ for the cuspidal curve $V^2+T^3=0$ in $\A^2_{(V,T)}$, 
which is non-normal. 
The quotient morphism $\pi : \tilde X \to Y$ is a trivial $\A^1$-bundle, but $Y\not\cong \A^3$.

At the present, whether $\tilde X$ is isomorphic to $\A^4$ or not is unknown.
If $\tilde X=X\times \A^1\cong \A^4$, we obtain a negative answer to the Cancellation Problem
which asks whether $X\times \A^1\cong \A^{n+1}$ implies $X\cong \A^n$ or not. 

\bigskip

Theorem \ref{Theorem 4.3} leads to the following question.

\medskip

\noindent{\bf Question}\; 
Let $X$ be a smooth acyclic affine fourfold with a free $G_a$-action.
Under the assumptions (1)--(3) in Theorem \ref{Theorem 4.3}, $X\cong \A^4$?

\medskip
\bigskip





\bibliographystyle{plain}
\bibliography{refs}

@article{BB,
  author		= "Bia{\l}ynicki-Birula, A.",
  title			= "On fixed point schemes of actions of multiplicative and additive groups",
  journal		= "Topology",
  volume		= "12",
  number		= "",
  pages			= "99--103",
  year			= "1973",
}

@article{BhD,
  author		= "Bhatwadekar, S. M. and Daigle, D.",
  title			= "On finite generation of kernels of locally nilpotent {$R$}-derivations of {$R[X,Y,Z]$}",
  journal		= "J. Algebra",
  volume		= "322",
  number		= "",
  pages			= "2915--2926",
  year			= "2009",
}

@article{Bon,
  author		= "Bonnet, P.",
  title			= "Surjectivity of quotient maps for algebraic ({$\C$},+)-actions",
  journal		= "Transform. Groups",
  volume		= "7",
  number		= "",
  pages			= "3--14",
  year			= "2002",
}

@article{DFN,
  author		= "Daigle, D. and Freudenburg, G. and Nagamine, T.",
  title			= "Generalizations of {S}amuel's criteria for a ring to be a unique factorization domain",
  journal		= "J. Algebra",
  volume		= "594",
  number		= "",
  pages			= "271--306",
  year			= "2022",
}

@article{DK,
  author		= "Daigle, D. and Kaliman, S.",
  title			= "A note on locally nilpotent derivations and variables of {$k[X,Y,Z]$}",
  journal		= "Canad. Math. Bull.",
  volume		= "52",
  number		= "4",
  pages			= "535--543",
  year			= "2009",
}

@article{ER,
  author		= "van den Essen, A. and van Rossum, P.",
  title			= "Coordinates in two variables over a {$\Q$}-algebra",
  journal		= "Trans. Amer. Math. Soc.",
  volume		= "356",
  number		= "5",
  pages			= "1691--1703",
  year			= "2004",
}

@article{Fuj,
  author		= "Fujita, T.",
  title			= "On the topology of non-complete algebraic surfaces",
  journal		= "J. Fac. Sci. Univ. Tokyo, Sec. IA",
  volume		= "29",
  number		= "",
  pages			= "503--566",
  year			= "1982",
}

@article{GMM,
  author		= "Gurjar, R.V. and Masuda, K. and Miyanishi, M.",
  title			= "{$\A^1$}-fibrations on affine threefolds",
  journal		= "J. Pure Appl. Algebra",
  volume		= "216",
  number		= "2",
  pages			= "296--313",
  year			= "2012",
}

@misc{GhoGupPal,
  author		= "Ghosh, Parnashree and Gupta, Neena and Pal, Ananya",
  title			= "On embedding of linear hyperplanes, ar{X}iv: 2405.07205v3",
  journal		= "",
  volume		= "",
  number		= "",
  pages			= "",
  year			= "",
}

@article{GMM2,
  author		= "Gurjar, R.V. and Masuda, K. and Miyanishi, M.",
  title			= "Affine threefolds with {$\A^2$}-fibrations",
  journal		= "Transform. Groups",
  volume		= "22",
  number		= "1",
  pages			= "187--205",
  year			= "2017",
}

@article{GKMMR,
  author		= "Gurjar, R.V. and Koras, M. and Masuda, K. and Miyanishi, M. and Russell,                           P",
  title			= "Affine threefolds admitting {$G_a$}-actions",
  journal		= "Math. Ann.",
  volume		= "373",
  number		= "3-4",
  pages			= "1211--1236",
  year			= "2017",
}

@incollection{GKMMR2,
  author		= "Gurjar, R.V. and Koras, M. and Masuda, K. and Miyanishi, M. and Russell,                           P",
  title			= "{$\A^1_*$}-fibrations on affine threefolds",
  editor		= "Masuda, K. and others",
  booktitle		= "Affine {A}lgebraic {G}eometry",
  pages			= "62--102",
  address		= "Singapore",
  publisher		= "World Sci. Publ.",
  year			= "2013",
}

@article{Kal,
  author		= "Kaliman, K.",
  title			= "Exotic analytic structures and {E}isenman intrinsic measures",
  journal		= "Israel Math. J.",
  volume		= "88",
  number		= "",
  pages			= "411--423",
  year			= "1994",
}

@article{Kal1,
  author		= "Kaliman, K.",
  title			= "Free {$\C_+$}-actions on {$\C^3$} are translations",
  journal		= "Invent. Math.",
  volume		= "156",
  number		= "",
  pages			= "163--173",
  year			= "2004",
}

@article{Kal2,
  author		= "Kaliman, S.",
  title			= "Polynomials with general {$\C^2$}-fibers are variables",
  journal		= "Pacific J. Math.",
  volume		= "203",
  number		= "",
  pages			= "161--190",
  year			= "2002",
}

@article{Kal3,
  author		= "Kaliman, S.",
  title			= "Proper {$G_a$}-actions on {$\C^4$} preserving a coordinate",
  journal		= "Algebra Number Theory",
  volume		= "12",
  number		= "2",
  pages			= "227--258",
  year			= "2018",
}

@article{KS,
  author		= "Kaliman, S. and Saveliev, N.",
  title			= "{$\C^+$}-actions on contractible threefolds",
  journal		= "Mich. Math. J.",
  volume		= "52",
  number		= "3",
  pages			= "619--625",
  year			= "2004",
}

@article{KZ,
  author		= "Kaliman, S. and Zaidenberg, M.",
  title			= "Affine modifications and affine hypersurfaces with a very transitive automorphism group",
  journal		= "Transform. Groups",
  volume		= "4",
  number		= "",
  pages			= "53--95",
  year			= "1999",
}

@article{Ma2,
  author		= "Masuda, K.",
  title			= "Factorial affine {$G_a$}-varieties with height one plinth ideals",
  journal		= "Transform. Groups",
  volume		= "27",
  number		= "4",
  pages			= "1287--1305",
  year			= "2023",
}

@incollection{Miy,
  author		= "Miyanishi, M.",
  title			= "Normal affine subalgebra of a polynomial ring",
  editor		= "",
  booktitle		= "Algebraic and topological theories (Kinosaki, 1984)",
  pages			= "37--51",
  address		= "Tokyo",
  publisher		= "Kinokuniya",
  year			= "1986",
}

@incollection{Miy4,
  author		= "Miyanishi, M.",
  title			= "Algebraic characterizations of the affine {$3$}-space",
  editor		= "",
  booktitle		= "Proceedings of the Algebraic Geometry Seminar, Singapore 1987",
  pages			= "53--67",
  address		= "Singapore",
  publisher		= "World Scientific",
  year			= "1988",
}

@article{Nag,
  author		= "Nagata, M.",
  title			= "A remark on the unique factorization theorem",
  journal		= "J. Math. Soc. Japan",
  volume		= "9",
  number		= "",
  pages			= "143--145",
  year			= "1957",
}

@article{Nor,
  author		= "Nori, M.",
  title			= "Zariski's conjecture and related problems",
  journal		= "Ann. Sci. Ecole Norm. Sup.",
  volume		= "16",
  number		= "2",
  pages			= "305--344",
  year			= "1983",
}

@article{Sat,
  author		= "Sathaye, A.",
  title			= "Polynomial ring in two variables over a {DVR}: a criterion",
  journal		= "Invent. Math.",
  volume		= "74",
  number		= "1",
  pages			= "159--168",
  year			= "1983",
}

@article{Win,
  author		= "Winkelmann, J.",
  title			= "On free holomorphic {$\C$}-actions on {$\C^n$} and homogeneous                                     {S}tein manifolds",
  journal		= "Math. Ann.",
  volume		= "286",
  number		= "",
  pages			= "593--612",
  year			= "1990",
}

@article{Zar,
  author		= "Zariski, O.",
  title			= "Interpr\'etations alg\'ebrico-g\'eom\'etriques du quatorzi\`eme probl\'eme de {H}ilbert",
  journal		= "Bull. Sci. Math. ",
  volume		= "78",
  number		= "2",
  pages			= "155--168",
  year			= "1954",
}

@book{Fre,
  author		= "Freudenburg, G.",
  title			= "Algebraic {T}heory of {L}ocally {N}ilpotent {D}erivations",
  address		= "Berlin Heidelberg",
  publisher		= "Springer",
  year			= "2017",
  series                = "Encyclopaedia of Mathematical Sciences vol. 136, 2nd edition",
}

@book{GMM3,
  author		= "Gurjar, R.V. and Masuda, K. and Miyanishi, M.",
  title			= "Affine {S}pace {F}ibrations",
  address		= "Berlin",
  publisher		= "De Gruyter",
  year			= "2021",
  series                = "De Gruyter Studies in Mathematics {\bf 79}",
}

@book{Matsu,
  author		= "Matsumura, H.",
  title			= "Commutative {R}ing {T}heory",
  address		= "Cambridge",
  publisher		= "Cambridge University Press",
  year			= "1989",
  series                = "Cambridge studies in advanced mathematics vol. 8",
}

@book{Miy3,
  author		= "Miyanishi, M.",
  title			= "Open {A}lgebraic {S}urfaces",
  address		= "Providence, RI",
  publisher		= "Amer. Math. Soc.",
  year			= "2001",
  series                = "CRM Monograph Series 12",
}

\end{document}